\documentclass[reqno,11pt]{amsart}
\usepackage[margin=2.5cm]{geometry}

\usepackage{amssymb}
\usepackage{amsmath}
\usepackage{amsfonts}
\usepackage{mathtools}
\usepackage{esint}
\usepackage{graphicx}
\usepackage{stackrel}
\usepackage{xcolor}
\usepackage{xypic}
\usepackage{lineno}

\usepackage{tikz}
\usepackage{tikz-cd}

\usepackage[shortlabels]{enumitem}
\usepackage{comment}
\usepackage{amsthm}
\usepackage{hyperref}
\hypersetup{
	colorlinks=true,
	linkcolor=blue,
	citecolor=blue,
	urlcolor=blue,
	pdfauthor={},
	pdftitle={}
}

\usepackage[capitalise, nameinlink, noabbrev]{cleveref} 

\newtheorem{thm}{Theorem}[section]

\newtheorem{lem}[thm]{Lemma}
\newtheorem{pro}[thm]{Proposition}
\newtheorem{theorem}{Theorem}

\newtheorem{thm.alpha}{Theorem}  
\newtheorem{cor.alpha}[thm.alpha]{Corollary}  
\newtheorem{question.alpha}[thm.alpha]{Question} 
\newtheorem{prop.alpha}[thm.alpha]{Proposition} 

\numberwithin{equation}{section}

\theoremstyle{definition}

\newtheorem{exa}[thm]{Example}

\newtheorem{rem}[thm]{Remark}
\newtheorem*{ack}{Acknowledgments}

\DeclareMathOperator{\Susp}{Susp}

\DeclareMathOperator{\seif}{Seif}
\DeclareMathOperator{\sing}{sing}

    \newcommand{\bpt}{B\!\left(\mathrm{pt}\right)}

\begin{document}


\title[Generalized Seifert Spaces]{Classification of generalized Seifert fiber spaces}

\author[F.~Galaz-García]{Fernando Galaz-García}
\address[Galaz-García]{Department of Mathematical Sciences, Durham University, United Kingdom}
\email{fernando.galaz-garcia@durham.ac.uk}

\author[J.~Núñez-Zimbrón]{Jesús Núñez-Zimbrón}
\address[Núñez-Zimbrón]{Facultad de Ciencias, UNAM, Mexico}
\email{nunez-zimbron@ciencias.unam.mx}


\begin{abstract}
We provide a symbolic classification of generalized Sei\-fert fiber spaces, which were introduced by Mitsuishi and Yamaguchi in the classification of collapsing Alexandrov $3$-spaces. Additionally, we show that the canonical double branched cover of a non-manifold generalized Seifert fiber space is a Seifert manifold and compute its symbolic invariants in terms of those of the original  space.  
\end{abstract}

\maketitle

\section{Main Results}

Generalized Seifert fiber spaces first appeared in the work of Mitsuishi and Yamaguchi on the classification of collapsing Alexandrov $3$-spaces \cite{MitsuishiYamaguchi2015} and generalize Seifert $3$-manifolds, sharing a similar structure: they decompose into fibers over a $2$-orbifold, with the fibers being either circles or closed intervals. As in the manifold case, small tubular neighborhoods of the circle fibers are standard fibered tori. In contrast, small neighborhoods of the interval fibers are homeomorphic to $\bpt$, a fundamental non-negatively curved Alexandrov $3$-space that serves as a ``singular'' fibered torus and can be written as the boundary connected sum of two cones over $P^2$ (see Section \ref{S:Preliminaries} for a detailed definition).

A closed (i.e., compact and without boundary) Alexandrov $3$-space is homeomorphic to either a $3$-manifold or to a non-orientable $3$-manifold with finitely many boundary components homeomorphic to $P^2$, the real projective plane, that have been capped off with cones over $P^2$ \cite{GalazGuijarro3dim}; the latter spaces have also appeared in the literature under the name \emph{singular $3$-manifolds} \cite{Quinn81,HAM,HAS}. Conversely, any $3$-manifold or singular $3$-manifold is homeomorphic to some Alexandrov space  \cite{GalazGuijarro3dim}. 

By the prime and Jaco--Shalen--Johannson (JSJ) decomposition theorems and Perelman's proof of Thurston's geometrization conjecture, every closed orientable $3$-manifold decomposes canonically into pieces whose interior is either Seifert fibered or hyperbolic (see \cite{Aschenbrenner.Friedl.Wilton} Theorem 1.6.1 (JSJ) together with Theorem 1.7.6 (geometrization), or \cite[Section 2.3]{Porti}). In particular, closed $\mathrm{Sol}$ manifolds have JSJ  decompositions consisting entirely of Seifert fibered pieces (they are graph manifolds), even though they are not globally Seifert fibered. The other six non-hyperbolic geometries 
($S^3$, $E^3$, $S^2\times\mathbb{R}$, $H^2\times\mathbb{R}$, $\widetilde{\mathrm{SL}_2(\mathbb{R})}$, and $\mathrm{Nil}$) occur precisely on Seifert fibered manifolds (see \cite[Theorem 5.3]{Scott}). In the Alexandrov setting, closed $3$-spaces also admit a geometric decomposition into geometric pieces  and any non-manifold geometric piece must be a $\mathbb{Z}_2$-quotient of an orientable geometric $3$-manifold (its orientable double branched cover) with geometry $S^3$, $E^3$, $H^3$, $S^2\times\mathbb{R}$, or $H^2\times\mathbb{R}$ (see \cite[Proof of Theorem 1.5 and Theorem 1.6]{GalazGuijarro3dim}). Hence, if such a piece is not hyperbolic, its double branched cover admits a Seifert fibration. By \cite{Meeks.Scott} and \cite{Dinkelbach.Leeb}, the $\mathbb{Z}_2$-action is \emph{standard}, i.e., it must preserve the geometric structure, and is conjugate to an isometry of a locally homogeneous metric. If the $\mathbb{Z}_2$-action preserves the corresponding Seifert fibration of the double branched cover, then  the quotient has a generalized Seifert fibration. This is the case, for example, when the geometry is $H^2\times\mathbb{R}$, since isometries preserve the foliation by the lines $\{x\}\times \mathbb{R}$ on the model space and hence preserve the induced Seifert fibration on any closed quotient (see \cite{Scott} or \cite[Theorem 2.2]{Meeks.Scott}; cf.\ also \cite{Mecchia.Seppi,Tollefson1978}). In the Euclidean case, in contrast, a standard action need not preserve a chosen Seifert fibration (see \cite{Meeks.Scott}).
Therefore, it is natural to consider a classification of generalized Seifert fiber spaces, as they constitute fundamental components of Alexandrov $3$-spaces and play a central role in their collapse theory \cite{MitsuishiYamaguchi2015,galaz.guijarro.nunez.zimbron,MitsuishiYamaguchi.Boundary}. Our first main result is such a classification, extending the classical one originally due to Seifert \cite{Seifert}.


\begin{theorem}
\label{T:Classification-generalized-Seifert}
Let $p\colon X\to B$ be a generalized Seifert fibration  with $X$ closed and connected. Then the generalized Seifert fiber space $X$ is determined up to fiber-preserving homeomorphism by the set of invariants
\[
\left\{ b,\varepsilon,g, \iota, \left\{ (\alpha_i, \beta_i) \right\}_{i=1}^n \right\}.
\]
\end{theorem}

For concision, the precise definitions of the symbols in Theorem~\ref{T:Classification-generalized-Seifert} are stated in Section \ref{S:Preliminaries}. 
Some of these invariants coincide with those appearing in the classification of local circle actions on Alexandrov $3$-spaces \cite{GalazZimbronLocalS1}. However, Theorem \ref{T:Classification-generalized-Seifert} is not a special case of the classification of local circle actions on Alexandrov $3$-spaces \cite[Theorem B]{GalazZimbronLocalS1}, as the canonical ``fibration'' of the space $\bpt$ does not arise from any local circle action (see the discussion after \cite[Lemma 6.1]{GalazZimbronLocalS1}). 

Our strategy to prove Theorem \ref{T:Classification-generalized-Seifert} involves excising small neighborhoods of the interval fibers and gluing appropriate ``blocks'' to obtain an associated space that does admit a local circle action, and use the corresponding classification. As part of the proof, for each set of invariants we construct an orbifold Riemannian metric which is locally invariant for the corresponding local circle action.

Every closed non-manifold Alexandrov $3$-space $X$ is the base of a two-fold branched cover $\pi\colon \tilde{X} \to X$ whose total space $\tilde{X}$ is a closed orientable manifold 
and whose branching set is the set of non-manifold points of $X$, i.e., points which have neighborhoods homeomorphic to a cone over $P^2$ \cite{GalazGuijarro3dim}. In particular, there is an
orientation-reversing (piecewise linear) involution $\tau \colon \tilde{X} \to \tilde{X}$ with only isolated fixed points such that $X$
is homeomorphic to the quotient $\tilde{X} /\tau$ \cite{GalazGuijarro3dim}. The manifold  $\tilde{X}$ is unique up to homeomorphism and we refer to it as the \emph{double branched cover} of $X$.   

In Theorem \ref{T:Double-cover-is-seifert}, we show that the canonical double branched cover of a non-manifold generalized Seifert fiber space $X$ is a Seifert manifold $\tilde{X}$ where the Seifert fibration commutes with the double branched cover construction. We also compute the symbolic invariants of $\tilde{X}$ in terms of those of the original space $X$, finding, roughly speaking, that the invariants of $X$ are doubled in $\tilde{X}$.  

\begin{theorem}
\label{T:Double-cover-is-seifert}
Let $p\colon X\to B$ be a generalized Seifert fibration where $X$ is closed, connected, and not a manifold, 
and let $\pi\colon\tilde{X} \to X$ be the canonical double branched cover of $X$. Then there exists a Seifert fibration $\tilde{p}\colon \tilde{X}\to \tilde{B}$ commuting with the double branched cover. Moreover, if $X$ is determined by the invariants
\[
\left\{ b,\varepsilon,g, \iota, \left\{ (\alpha_i, \beta_i) \right\}_{i=1}^n \right\},
\]
then $\tilde{X}$ is determined by the invariants 
\[
\left\{\tilde{b} , \tilde{\varepsilon}, \tilde{g}, 0, \left\{ (\tilde{\alpha}_i, \tilde{\beta}_i) \right\}_{i=1}^{2n} \right\},
\]
where $\tilde{\varepsilon}=\varepsilon$,  $(\tilde{\alpha}_{2k},\tilde{\beta}_{2k})=(\alpha_k, \beta_k)$ and $(\tilde{\alpha}_{2k-1},\tilde{\beta}_{2k-1})=(\alpha_k, \beta_k)$ for each $k=1,2,\ldots,n$, and the base genera satisfy
\[
\chi(\tilde B)=2\chi(B)-\iota.
\]
In particular, if $B$ is closed and orientable, then $\tilde g=2g+\frac{\iota}{2}-1$ (hence $\iota$ is even). Moreover, the rational Euler number of the Seifert bundle $\tilde p\colon \tilde{X}\to \tilde{B}$ is $e(\tilde X)=0$. Equivalently, 
\[
\tilde b = -2\sum_{i=1}^n \frac{\beta_i}{\alpha_i}.
\]
In particular, if $n=0$, then $\tilde b=0$.
\end{theorem}

Note that the relations in Theorem~\ref{T:Double-cover-is-seifert} yield constraints for a given Seifert $3$-manifold to be the double branched cover of the total space of a generalized Seifert fibration.
 When we couple Theorem \ref{T:Double-cover-is-seifert} with the work of Cheeger and Gromov on collapsing Riemannian manifolds with bounded curvature \cite{cheeger.gromov.1986}, we conclude that any generalized Seifert fiber space collapses with bounded diameter and a uniform lower bound on the curvature. Moreover, together with the classification of Mitsuishi and Yamaguchi of closed collapsing Alexandrov $3$-spaces \cite{MitsuishiYamaguchi2015}, this characterizes Alexandrov $3$-spaces that collapse with bounded diameter and a uniform lower curvature bound. 
 
 Further, we observe that the Euler characteristic $\chi(B)$ together with $\iota$, $\tilde{b}$, and the remaining invariants in Theorem~\ref{T:Double-cover-is-seifert}, may be used to recover the topology of $\tilde{X}$ and, in concrete cases, also that of its quotient $X = \tilde{X}/\tau$  (compare with Example~\ref{ex:second.fibration}).
 We note that Bonahon and Siebenmann classified Seifert fibered $3$-orbifolds in full generality \cite{Bonahon.Siebenmann}. We work in the topological and Alexandrov settings. Our arguments are independent of the orbifold classification and rely instead on the classification of local circle actions \cite{FintushelLocal,GalazZimbronLocalS1,OrlikRaymondLocal}.
\hfill
\\

Our article is organized as follows. In Section~\ref{S:Preliminaries}, we recall basic facts on generalized Seifert fibered spaces and their topology, and collect basic facts on the classification of Alexandrov $3$-spaces with local circle actions that we will use in the proofs of the main theorems. We prove Theorems~\ref{T:Classification-generalized-Seifert} and \ref{T:Double-cover-is-seifert} in Sections~\ref{S:proof.theorem.A} and \ref{S:proof.theorem.B}, respectively.

\begin{ack}
JNZ acknowledges support from PAPIIT--UNAM projects IN101322 and IA103925. JNZ wishes to thank Diego Corro and Luis Jorge Sánchez Saldaña for useful conversations. Both authors thank the anonymous referee for a thorough review and helpful observations and suggestions.
\end{ack}

\section{Preliminaries}
\label{S:Preliminaries}

To keep our presentation concise, throughout the article we assume the reader to be familiar with the basic theory of Seifert manifolds and Alexandrov spaces of curvature bounded below, particularly in the three-dimensional case. For an introduction to  Alexandrov spaces, we refer to the standard references \cite{AlexanderKapovitchPetrunin, BuragoBuragoIvanov}. For basic results in dimension three, we refer the reader to \cite{GalazGuijarro3dim,FGGJNZ.Survey}. For further results on Alexandrov $3$-spaces, see  \cite{barcenas.nunez.zimbron.2021,barcenas.sedano.2023,Deng.et.al.3D.Alex.Ricci,FrancoGalazLarranagaGuijarroHeil,galaz.guijarro.nunez.zimbron,GalazZimbronLocalS1,galaz-garcia.searle.2011,galaz-garcia.tuschmann.2019,Kapovich,MitsuishiYamaguchi.Boundary,nunez.zimbron.2018.S1.actions}.  Our main references for the basic theory of Seifert fiber manifolds will be \cite{JankinsNeumann1983,OrlikBook,Scott}.
In this section, we only recall the key aspects necessary for our discussion. We will use the symbol ``$\cong$'' to denote homeomorphism between topological spaces.

\subsection{Alexandrov \texorpdfstring{$3$}{3}-spaces}
\label{ss:alexandrov.3-spaces}
Let $X$ be a closed Alexandrov $3$-space. The space of directions $\Sigma_x X$ at each point $x$ in $X$ is a closed Alexandrov $2$-space with curvature bounded below by $1$ and, by the Bonnet--Myers Theorem \cite[Theorem 10.4.1]{BuragoBuragoIvanov}, the fundamental group of $\Sigma_x X$ is finite. Thus, $\Sigma_x X$ is homeomorphic to a $2$-sphere $S^2$ or a real projective plane $P^2$. A point in $X$ whose space of directions is homeomorphic to $S^2$ is a \textit{topologically regular point}. Points in $X$ whose space of directions is homeomorphic to $P^2$ are \textit{topologically singular}.  We will denote the set of topologically singular points of $X$ by $\sing(X)$. The set of topologically regular points is open and dense in $X$. By Perelman's Conical Neighborhood Theorem \cite{Perelman}, every point $x$ in $X$ has a neighborhood that is pointed-homeomorphic to the cone over $\Sigma_x X$. Therefore, the set of topologically singular points of $X$ is finite, and $X$ is homeomorphic to a compact $3$-manifold with a finite number of $P^2$-boundary components to which one glues in cones over $P^2$. If $X$ contains topologically singular points, it is homeomorphic to the quotient of a closed, orientable, topological $3$-manifold $\tilde{X}$ by an orientation-reversing involution $\tau\colon\tilde{X}\to \tilde{X}$ with only isolated fixed points. The $3$-manifold $\tilde{X}$ is the \emph{orientable double branched cover of $X$} (see, for example, \cite[Lemma 1.7]{GalazGuijarro3dim}). One may lift the metric on $X$ to $\tilde{X}$ so that $\tilde{X}$ becomes an Alexandrov space with the same lower curvature bound as $X$ and $\tau$ is an isometry with respect to the lifted metric (see \cite[Lemma 1.8]{GalazGuijarro3dim} and \cite[Proposition 3.4]{Deng.et.al.3D.Alex.Ricci}, and compare with \cite[Lemma 5.2]{Grove.Wilking}). In particular, $\tau$ is equivalent to a smooth involution on $\tilde{X}$ considered as a smooth $3$-manifold.


\subsection{Generalized Seifert fiber spaces}
\label{ss:generalized.seifert.fiber.spaces}
Let us recall the definition of generalized Seifert spaces, introduced by Mitsuishi and Yamaguchi in \cite{MitsuishiYamaguchi2015}.
We begin with the definition of the space $\bpt$ (cf. \cite[Example 1.2]{MitsuishiYamaguchi2015}). Consider the isometric involution (with respect to the flat metric) 
\begin{eqnarray}
\label{involution-solid-torus-to-bpt}
\alpha\colon S^1\times D^2 &\to& S^1\times D^2\\
(e^{i\theta},x)&\mapsto& (e^{-i\theta},-x).
\end{eqnarray}
The space $\bpt$ is then defined as $\left(S^1\times D^2\right)/\left\langle \alpha \right\rangle$, and we denote the projection by $\pi\colon S^1\times D^2 \to \bpt$. There is a natural projection 
\[
q\colon B(\mathrm{pt}) \to K_1(S^1(1/2)) 
\]
induced by $(e^{i\theta},x)\mapsto x$. Here, $K_1(S^1(1/2))$ is the ball of radius one centered at the vertex of a Euclidean cone
over $S^1(1/2)$, a circle in $\mathbb{R}^2$ of radius $1/2$ or, equivalently, a circle of length $\pi$.


\begin{rem}
    Observe that the restriction of $\alpha$ to $S^1\times S^1$ has as quotient the flat Klein bottle $K^2$, viewed as the non-orientable circle bundle over $S^1$ (see, for example, \cite[Lemma 3]{Natsheh}). Consequently, the fibration on $K^2$ induced by the restriction $q\colon K^2 \to S^1$ is fiberwise equivalent (i.e., there is a fiberwise homeomorphism) to the fibration induced by the free local circle action (see next section for the definition of a local circle action) on $K^2$ given by rotating one of the circle factors.  
\end{rem}

A \textit{generalized Seifert fibration} of a (topological) $3$-orbifold $X$ over a connected (topological) $2$-orbifold $B$ (possibly with non-empty or disconnected boundary) is a map $f\colon X\to B$ whose fibers are homeomorphic to circles or bounded closed intervals, and which satisfies the following conditions. For every $x\in B$, there is a neighborhood $U_x$ homeomorphic to a $2$-disk satisfying the following conditions:
\begin{itemize}
    \item[(i)] If $f^{-1}(x)$ is homeomorphic to a circle, then there is a fiber-preserving homeomorphism of $f^{-1}(U_x)$ to a Seifert fibered solid torus in the usual sense (see \cite[Section 1.7]{OrlikBook} for the definition of a Seifert fibered solid torus and the associated \textit{Seifert invariants} $(\alpha,\beta)$).
    In this case, we say that $f^{-1}(x)$ is a \emph{$C$-fiber}, and that $x$ \emph{has $C$-fiber type}. Observe that, in a sufficiently small tubular neighborhood of the fiber, the Seifert fibration is given by the orbits of an effective circle action which is determined by the invariants $(\alpha, \beta)$. Slightly abusing terminology, we define the \emph{isotropy} of the central $C$-fiber as the isotropy group of the fiber considered as an orbit of the effective circle action inducing the corresponding Seifert fibration on the Seifert fibered solid torus. 
     
    \item[(ii)] If $f^{-1}(x)$ is homeomorphic to an interval, then there exists a fiber-preserving homeomorphism of $f^{-1}(U_x)$ to the space 
    $B(\mathrm{pt})$, with respect to the fibration
    \[
    \left(B(\mathrm{pt}), q^{-1}(o)\right) \to \left( K_1(S^1(1/2)), o\right).
    \]
    In this case, we say that $f^{-1}(x)$ is an \emph{$I$-fiber}, and that  $x$ \emph{has $I$-fiber type}.
\end{itemize}
Moreover, for any compact component $C$ of the boundary of $B$, we require the existence of a collar neighborhood $N$ of $C$ in $B$ such that $f|_{f^{-1}(N)}$ is a circle bundle over $N$. In this context, we say that $X$ is a \emph{generalized Seifert fiber space} with base $B$ and use the notation $X=\seif(B)$. The base spaces of two generalized Seifert fibered spaces $X=\seif(B)$ and $Y=\seif(C)$ are \emph{isomorphic} if there exists a homeomorphism between $B$ and $C$ that preserves the fiber type.

As noted in \cite[Section 6]{GalazZimbronLocalS1}, when there are $I$-fibers, the fibration $f$ is not induced by any local circle action. Note that when $X$ is compact and has no boundary, the base $B$ is homeomorphic to a closed surface.

\subsection{Local circle actions}

We briefly summarize the classification of local circle actions on Alexandrov $3$-spaces obtained in \cite{GalazZimbronLocalS1}, which will play a key role in our discussion (see \cite{FintushelLocal, OrlikRaymondLocal} for the manifold case). 

A \emph{local circle action on a closed Alexandrov $3$-space $X$} is a decomposition of $X$ into (possibly one-point) disjoint, simple, closed curves (the \textit{fibers}), having a tubular neighborhood which admits an effective circle action whose orbits are the curves of the decomposition. 
The local  circle action is \emph{by isometries} if the circle actions on each tubular neighborhood of the fibers are by isometries with respect to the restricted metric.


\subsection*{Fiber types} We define the fiber types of a local circle action by the isotropy of the local action, as well as the isotropy representation and if any point on the orbit is a regular or singular point (see \cite[Section 3]{GalazZimbronLocalS1}). The possible types are the following:

\begin{itemize}
    \item \emph{$F$-fibers} are topologically regular fixed-point fibers.
    \item \emph{$SF$-fibers} are topologically singular fixed-point fibers.
    \item \emph{$E$-fibers} are those that correspond to $\mathbb{Z}_k$ isotropy, acting in such a way that local orientation is preserved.
    \item \emph{$SE$-fibers} correspond to $\mathbb{Z}_2$ isotropy, reversing the local orientation.
    \item \emph{$R$-fibers} are fibers that are not $F$-, $SF$-, $E$- or $SE$-fibers.
\end{itemize}  

The sets of points in $F$-, $SF$-, $E$-, $SE$- and $R$-fibers are denoted by $F$, $SF$, $E$, $SE$ and $R$, respectively. The quotient space $X^{\ast}$ of the local circle action, with the quotient topology, admits a $2$-orbifold structure (see
for example \cite{GalazZimbronLocalS1}). Its boundary consists of the images of $F$-, $SF$- and $SE$-fibers under the natural projection map, while the interior consists of $R$-fibers and a finite number of $E$-fibers.

\subsection*{Building blocks} A closed Alexandrov $3$-space $X$ with a local isometric circle action can be decomposed into \emph{building blocks} of types $F$, $SF$, $E$ and $SE$, defined by considering small invariant tubular neighborhoods of connected components of the strata in $X$ consisting of fibers of the corresponding type. A building block is \emph{simple} if its boundary is homeomorphic to a torus, and \textit{twisted} if its boundary is homeomorphic to a Klein bottle. Note that the stratum of $R$-fibers is an $S^1$-fiber bundle with structure group $\mathrm{O}(2)$.

\begin{thm}[cf.\ \protect{\cite[Theorem B]{GalazZimbronLocalS1}}]
\label{THM:INVARIANTS}
Let $X$ be a closed, connected Alexandrov $3$-space with $2r\geq 0$ topologically singular points. Then, the set of isometric local circle actions (up to equivariant equivalence) is in one-to-one correspondence with the set of unordered tuples
\[
\left\{b; \varepsilon, g, (f,k_1), (t,k_2), (s,k_3); \{ (\alpha_i, \beta_i) \}_{i=1}^n; (r_1,r_2, \ldots, r_{s-k_{3}}); (q_1, q_2, \ldots, q_{k_3})\right\}.
\]
\end{thm}

The definition of the symbols appearing in this result is as follows:
\begin{itemize}
     \item The numbers $k$, $k_1$, $k_2$, and $k_3$ satisfy $k_1 + k_2 + k_3 = k$.
    \item $(\varepsilon,k)$ is the pair that classifies the $S^1$-bundle of $R$-fibers according to \cite[Theorem 3.2]{GalazZimbronLocalS1}.
    \item $g\geq 0$ is the genus of $X^{*}$.
    \item The symbols $f, t, k_1, k_2$ are non-negative integers  such that $k_1 \leq f$ and $k_2\leq t$, where $k_1$ is the number of twisted $F$-blocks and $k_2$ is the number of twisted $SE$-blocks (therefore $f-k_1$ is the number of simple $F$-blocks and $t-k_2$ is the number of simple $SE$-blocks).
    \item $n$ is the number of $E$-fibers and we let $\{ (\alpha_i, \beta_i)\}_{i=1}^n$ be the corresponding Seifert invariants (see \cite{OrlikRaymondLocal} for the definition).

    \item We let $s, k_3$ be non-negative integers, where  $k_3\leq s$ is the number of twisted $SF$-blocks (thus $s-k_3$ is the number of simple $SF$-blocks).
    \item $(r_1, r_2, \ldots, r_{s-k_{3}})$ and $(q_1, q_2, \ldots, q_{k_3} )$ are $(s-k_{3})$- and $k_3$-tuples of non-negative even integers corresponding to the number of topologically singular points in each  simple and twisted $SF$-block, respectively.

     \item The invariant $b$ is an integer (or integer mod $2$) representing an obstruction class in the following way: Let $X_0$ denote the subset of $X$ obtained by removing from $X$ sufficiently small (open) tubular neighborhoods of the $F$-, $SE$-, $E$-, and $SF$-fibers which are disjoint. Then $X_0$ consists of $R$-fibers. Observe that the Seifert invariants determine cross-sections on the boundaries of $X_0$ which correspond to tubular neighborhoods of $E$-fibers and that any section that could already be defined on $X_0$ can be extended to the tubular neighborhoods of $F$-, $SF$-, and $SE$- fibers (see \cite[Proof of Theorem B]{GalazZimbronLocalS1}). Thus we define $b$ as the obstruction class in 
     \[
     H^2(X^{\ast} \setminus \mathrm{int}(V_1^{\ast}\cup\ldots\cup V_n^{\ast}), \partial (V_1^{\ast}\cup\ldots\cup V_n^{\ast}) ;\mathbb{Z}),
     \]
         where the $V_i$ are the tubular neighborhoods of the $E$-orbits, to extending a cross-section defined on $\partial (V_1^{\ast}\cup\ldots\cup V_n^{\ast})$ to all $X_0$. It follows from \cite{GalazZimbronLocalS1} and the manifold case in \cite{OrlikRaymondLocal}, \cite{FintushelLocal} that $b$ is an integer (or an integer mod $2$ depending on the orientability of $X^{\ast}$), subject to the following restrictions: $b=0$ if $f+t>0$ or if $\varepsilon\in\{o_2,n_1,n_3,n_4\}$ and some $\alpha_i=2$ (see \cite[Theorem 3.2]{GalazZimbronLocalS1} for the precise definitions of the $o_i$, $n_j$ and $\alpha_l$); $b\in\{0,1\}$ if $f+t=0$ and $\varepsilon\in\{ o_2,n_1,n_3,n_4\}$ and all $\alpha_i\neq 2$. In the remaining cases, $b$ is an arbitrary integer. It is worth pointing out that $b$ plays a role, for example, when  $SF\cup F\cup SE \cup E=\emptyset$ and that the structure group reduces to $\mathrm{SO}(2)$ (see, for example, the descriptions of $b$ in Lemmas 2 and 3 in \cite[Section 1.9]{OrlikBook}).
   
\end{itemize}

\section{Invariants and Classification}
\label{S:proof.theorem.A}

In this section, we define the equivariant and topological symbols that appear in Theorem \ref{T:Classification-generalized-Seifert} and prove this result.

\subsection{Invariants}
Let $X=\seif(B)$ be a closed, generalized Seifert fiber space and associate the following information to $X$:\\

\begin{itemize}
    \item The subset $X_0$ of $X$ obtained by removing all $I$-fibers and $C$-fibers with nontrivial isotropy is an $S^1$-bundle with structure group $\mathrm{O}(2)$, so we associate to $X$ the classifying pair $(\varepsilon,k)$. In fact, it will follow from the proof of Theorem \ref{T:Classification-generalized-Seifert} that $k=\iota$ (where $\iota$ is defined below). 
     We point out that, if instead we remove small, disjoint invariant neighborhoods of the $I$- and $C$-fibers with nontrivial isotropy, the resulting set is a deformation retract of $X_0$ and the classifying invariants $(\varepsilon,k)$ do not change. In what follows, we will also denote this deformation retract by $X_0$.
    \item The class $b$ is defined in a similar way to that of the local circle action case:  We remove small, pairwise disjoint tubular neighborhoods of the $I$-fibers and those $C$-fibers with nontrivial isotropy, producing a subset $X_0\subset X$. Then, the Seifert invariants of the $C$-fibers with finite, nontrivial isotropy specify cross-sections to the fibration on the boundaries of the tubular neighborhoods. Thus, $b$ is defined as the obstruction class to extending these cross-sections to a full cross-section (compare with Lemma~\ref{L:double.branched.cover} below). As a consequence, any cross-section $s_\partial\colon \partial\bpt^* \to \partial \bpt$ of the boundary fibration $\partial \bpt \to \partial \bpt^*$ (which is the nontrivial circle bundle over $S^1=\partial\bpt^*$) can be extended to a cross-section $s\colon \bpt^*\to \bpt$ of the fibration $\bpt \to \bpt^*$. 
    
    \item The symbol $g$ stands for the genus of $B$.
    \item The symbol $\iota$ is the (necessarily finite) number of $I$-fibers.
    \item Finally, the ordered pairs $(\alpha_i,\beta_i)$ are the Seifert invariants of the $C$-fibers with nontrivial isotropy.
\end{itemize}

All of these symbols, except $\iota$, are inspired by those assigned to a space with a local circle action as can readily be seen. In turn, they are subject to the same restrictions on its values as those indicated in Theorem \ref{THM:INVARIANTS} (except obviously $\iota$, which is new); however, not all values are possible under the assumption that $X$ is a generalized Seifert space, 
as will follow from the proofs of Theorem \ref{T:Classification-generalized-Seifert} and \ref{T:Double-cover-is-seifert}. For example, if $B$ is closed and orientable, $\iota$ must be even. 
 
To summarize, to each generalized Seifert bundle we associate a set of invariants of the form 
\begin{equation}
\label{list-of-invariants}
\left\{ b,(\varepsilon,k),g, \iota, \left\{ (\alpha_i, \beta_i) \right\}_{i=1}^n \right\}.
\end{equation}

\subsection{Classification}

If there are no $I$-fibers, Theorem \ref{T:Classification-generalized-Seifert} reduces to the classification of usual Seifert manifolds by letting $\iota=0$, and the fact that the classical classification of Seifert (see \cite{Seifert}) is contained in the classification of local circle actions on $3$-manifolds of Fintushel, Orlik and Raymond (see \cite{FintushelLocal, OrlikRaymondLocal}). Therefore, to prove Theorem \ref{T:Classification-generalized-Seifert},  it is sufficient to consider the case in which the set of $I$-fibers is non-empty. Note that this implies that $\sing(X)$, the set of topologically singular points of $X$, is non-empty.\\

To establish Theorem \ref{T:Classification-generalized-Seifert}, we need to prove existence and uniqueness statements as follows: that given a list of symbols as in \ref{list-of-invariants}, one is able to construct a generalized Seifert fiber space realizing these symbols, and that if two generalized Seifert fiber spaces share the same set of invariants, then there exists a fiber-preserving homeomorphism  between them.

\subsection*{Proof of Theorem \ref{T:Classification-generalized-Seifert}} 
We first prove the uniqueness statement. Let $X_1=\seif(B_1)$ and $X_2=\seif(B_2)$ be closed and connected generalized Seifert spaces and let $p_j\colon X_j\to B_j$, $j=1,2$ denote the generalized Seifert fibrations. It is clear that if $X_1$ is fiberwise homeomorphic to $X_2$, then $B_1$ is isomorphic to $B_2$. We now prove that if there is an isomorphism of fiber spaces $\varphi\colon B_1 \to B_2$, then $X_1$ is equivalent to $X_2$.\\

Let $I^{\ast}_{j}=\{x^{\ast}_{j1}, x^{\ast}_{j2}, \ldots, x^{\ast}_{j\iota}\}\subseteq B_j$ for $j=1,2$, be the collection of points of $I$-fiber type on $B_j$. We now choose disjoint closed disk neighborhoods $D_{1i}$ of each $x^{\ast}_{1i}$ such that $p_1^{-1}(D_{1i})\cong \bpt$. Without loss of generality, as $\varphi$ is fiber-type preserving, we can assume $\varphi(x^{\ast}_{1i})=x^{\ast}_{2i}$ and denote $D_{2i}= \varphi(D_{1i})$, $i=1,\ldots, \iota$. Furthermore, we can assume that $D_{1i}$ is small enough so that $p_2^{-1}(D_{2i})\cong \bpt$ for all $i=1,\ldots,\iota$. 

Now, observe that $\varphi$ restricts to a fiber-type preserving homeomorphism $B_1 \setminus \bigcup_{i=1}^{\iota} \mathrm{int}(D_{1i}) \to B_2 \setminus \bigcup_{i=1}^{\iota} \mathrm{int}(D_{2i})$ and consider for $j=1,2$
\[
Z_j:= X_j\setminus \bigcup_{i=1}^{\iota} p_j^{-1}\left( \mathrm{int}(D_{ji}) \right). 
\]
Then $\partial Z_j$ is a disjoint union of $\iota$ copies of Klein bottles $K^2$. Observe moreover that $Z_j$ is a Seifert manifold with boundary. 

We let $V_{ji}$ be twisted $F$-blocks for $j=1,2$, $i=1,\ldots, \iota$ and consider the spaces $Y_j$ obtained by gluing each $V_{ji}$ to each boundary component of $Z_j$ via fiber-preserving homeomorphisms. To this end, recall from \cite{GalazZimbronLocalS1} (see the subsection ``Block Types'') that the fibrations restricted to $\partial V_{ji}$ are unique up to fiber-preserving homeomorphism; note that this is also the case for the fibration on the boundary of $\bpt$. Further, note that the spaces $Y_j$ satisfy that $\sing(Y_j)=\emptyset$ as all topologically singular points of $X_j$ must be on $I$-fibers and not on $C$-fibers. Thus, each $Y_j$ is a closed topological $3$-manifold, and carries a local circle action. Therefore, by the work of Orlik and Raymond \cite{OrlikRaymondLocal} and Fintushel \cite{FintushelLocal}, the $Y_j$ are determined up to fiber-preserving homeomorphism by the same set of invariants, namely 
\[
\left\{b,\varepsilon,g, (f=\iota,k_1=\iota), (t=0,k_2=0), \left\{(\alpha_j,\beta_j) \right\}_{j=1}^{n} \right\},
\]
and there is a fiber-preserving homeomorphism $\Psi\colon Y_1 \to Y_2$ which is determined by cross-sections to the local actions away from the exceptional orbits. 

Observe now that $\Psi(V_{1i})\subset Y_2$ is a twisted $F$-block so that, without loss of generality, we can assume $\Psi(V_{1i})=V_{2i}$, for each $i=1,\ldots,\iota$. Thus, by restricting $\Psi$ we obtain a fiber-preserving homeomorphism $\Psi|_{Z_1}\colon Z_1 \to Z_2$ determined by cross-sections $\rho_j\colon \partial ( B_j \setminus \bigcup_{i=1}^{\iota} \mathrm{int}(D_{ji}))\to \partial Z_j$, that is, such that 
\[
\Psi\left( \rho_1\left( \partial ( B_1 \setminus \bigcup_{i=1}^{\iota} \mathrm{int}(D_{1i}))\right)\right) = \rho_2\left( \partial ( B_2 \setminus \bigcup_{i=1}^{\iota} \mathrm{int}(D_{2i}))\right).
\]
Furthermore, as the local circle action restricted to $\partial Z_j$ is equivalent to the ``free'' local circle action, we can parametrize each component of $\partial Z_j$, homeomorphic to the Klein bottle, seen as the nontrivial bundle $S^1\widetilde{\times}S^1$ as an equivalence class (with respect to the involution $\alpha$ as in \eqref{involution-solid-torus-to-bpt}) of points $[e^{i\theta}, e^{i\eta}]$, where $e^{i\theta}\in \rho_1\left( \partial ( B_1 \setminus \bigcup_{i=1}^{\iota} \mathrm{int}(D_{1i}))\right) \cong S^1$ and $e^{i\eta}$ parametrizes the fiber direction. Hence, $\Psi|_{\partial Z_1}$ is of the form 
\[
\Psi|_{\partial Z_1}([e^{i\theta}, e^{i\eta}]) = [\Psi_1(e^{i\theta}), \Psi_{2,\theta}(e^{i\eta})],
\]
where $\Psi_1$ and $\Psi_{2,\theta}$ are self-homeomorphisms of $S^1$ for each $\theta \in [0,2\pi)$. 

We now conclude by observing that $\Psi$ extends to a fiber-preserving homeomorphism $X_1\to X_2$ by parametrizing each $\bpt$ block by classes of points  $[e^{i\theta}, re^{i\eta}]$ where $r\in [0,1]$, and defining 
\[
\Psi|_{\bpt}([e^{i\theta}, re^{i\eta}]):= [\Psi_1(e^{i\theta}), r\Psi_{2,\theta}(e^{i\eta})].
\]

We now prove the existence part of the theorem. By the classification of spaces with local circle actions \cite{GalazZimbronLocalS1}, there exists a space $Y$ with a local circle action with the invariants
\[
\left\{b,\varepsilon,g, (f=\iota,k_1=\iota), (t=0,k_2=0), \left\{(\alpha_j,\beta_j) \right\}_{j=1}^{n} \right\},
\]
and an orbifold Riemannian metric of curvature bounded below. 
In a similar fashion to the uniqueness part of the proof, we take out all $\iota$ twisted $F$-blocks from $Y$, obtaining a Riemannian orbifold with boundary $Z$, curvature bounded from below, and a local circle action. We now take $\iota$ copies $V_i$ of $\bpt$ with the orbifold metric induced by the flat metric on $S^1\times D^2$. 

 By the smooth Collar Neighborhood Theorem \cite[Theorem 9.25]{Lee}, and since the boundary components of $Z$ and each $V_i$ are Klein bottles contained in their smooth loci,
$\partial Z$ and $\partial V_i$ admit smooth product collars, and we may assume the metrics to be product near the boundary and that their restrictions to corresponding components of $\partial Z$ and $\partial V_i$ are isometric.
Therefore, we can glue the $V_i$ to $Z$ fiber-preservingly by $\iota$ copies of the cylinder $I\times K^2$ in such a way that the metric on $\{0\}\times K^2$ coincides with the metric on the $i$-th boundary component of $Z$, and the metric on $\{1\}\times K^2$ coincides with the metric on the boundary of $V_i$. Moreover, this can be done so that the glued space $X$ is a smooth Riemannian orbifold which, by compactness, has curvature bounded below, and in particular, $X$ is an Alexandrov space admitting a generalized Seifert fiber space structure with the invariants
\[
\left\{ b,(\varepsilon,k),g, \iota, \left\{ (\alpha_i, \beta_i) \right\}_{i=1}^n \right\}.
\]
As we pointed out, $k=\iota$ so that we remove the superfluous invariant $k$ and uniquely associate (up to fiber-preserving homeomorphism) the list 
\begin{equation}
\left\{ b,\varepsilon,g, \iota, \left\{ (\alpha_i, \beta_i) \right\}_{i=1}^n \right\}.
\end{equation}
to each generalized Seifert fiber space, as stated in Theorem \ref{T:Classification-generalized-Seifert}.
\qed
\section{Compatibility with the double branched cover}
\label{S:proof.theorem.B}
The canonical double branched cover of an Alexandrov space $X$ with $\sing(X)\neq \emptyset$ is a fundamental construction which, in dimension $3$, de-singularizes the space \cite{GalazGuijarro3dim, HarveySearle}. In this section, we show that the double branched cover of a generalized Seifert fiber space is a Seifert manifold and compute its invariants from the invariants of the original space, thus proving Theorem~\ref{T:Double-cover-is-seifert}.\\ 

The following lemma follows directly from the fact that the orientation double cover of a non-orientable manifold can be characterized in the following way: Let $M$ be a non-orientable manifold and $p\colon\tilde{M}\to M$ its orientation double cover. If $N$ is an oriented manifold and $\pi\colon N\to M$ is a double cover with an orientation-reversing nontrivial deck transformation, then $(N,\pi)$ is isomorphic to $(\tilde{M},p)$, (see for example \cite{Kreck}). 
 Recall the canonical construction of the orientable double branched cover $\tilde{X}$ of a topologically singular Alexandrov $3$-space $X$: if $X$ is a closed Alexandrov $3$-space whose set of topologically singular points is $\sing(X)\neq \varnothing$, then let $U_1,\ldots,U_k$ be pairwise disjoint open cone neighborhoods of the points in $\sing(X)$, and set
$X_0 = X\setminus \bigcup_{i=1}^{k}U_i$. Then $X_0$ is a non-orientable $3$-manifold whose boundary consists of finitely many copies of $P^2$, and we obtain the canonical orientable double branched cover $\pi\colon \tilde{X}\to X$ by taking the orientation double cover $\tilde{X}_0$ of $X_0$ and then capping off each lifted $S^2$-boundary component  by a $3$-ball, so the deck involution extends to $\tilde{X}$ with isolated fixed points projecting to points in $\sing(X)$ (cf.\ \cite{GalazGuijarro3dim}). 

\begin{lem}
\label{L:double.branched.cover}
The (canonical) orientable double branched cover of $\bpt$ is isomorphic, as a double branched cover (i.e., via a $\mathbb{Z}_2$-equivariant homeomorphism commuting with the covering projections), to $S^1\times D^2$ with deck involution $\alpha\colon (e^{i\theta}, x)\mapsto (e^{-i\theta}, -x)$, and the branched covering map $\pi$ is the quotient map $S^1\times D^2 \to (S^1\times D^2)/\langle \alpha \rangle=\bpt$.
\end{lem}

 Note that the involution $\alpha$ has exactly two fixed points $(\pm 1,0)$. Their images in $\bpt$ are the two branch points of $\pi$, and they are the endpoints of the $I$-fiber $q^{-1}(o)$ in the generalized Seifert fibration $q\colon \left(\bpt, q^{-1}(o)\right) \to \left( K_1(S^1(1/2)), o\right)$.

With this observation in hand, we may prove our second main theorem.

\subsection*{Proof of Theorem~\ref{T:Double-cover-is-seifert}} 
We divide the proof in two parts. First, we show the existence of the lifted Seifert fibration $\widetilde{p}\colon \widetilde{X}\to \widetilde{B}$. We then determine the Seifert invariants of this fibration.


\subsubsection*{The lifted Seifert fibration.} The existence and properties of the lifted Seifert fibration are given in the following proposition.


\begin{pro}
\label{prop:lift_to_orientable_branched_cover}
Let $p\colon X \to B$ be a generalized Seifert fibration with $X$ closed, connected, and not a manifold. Let $\sing(X) \subset X$ denote the set of topologically singular points of $X$  and let
\[
\pi\colon \widetilde X \longrightarrow X
\]
be the canonical orientable double branched cover, branched precisely along $\sing(X)$. Then there exists a compact $2$-orbifold $\widetilde B$, a double branched cover
\[
\widetilde\pi\colon \widetilde B \longrightarrow B,
\]
and a Seifert fibration
\[
\widetilde p\colon \widetilde X \longrightarrow \widetilde B,
\]
such that:
\begin{enumerate}
\item $\widetilde\pi$ is branched precisely over $p(\sing(X)) \subset B$.
\item The following diagram commutes:
\begin{equation}
\begin{tikzcd}
\widetilde X \arrow[r,"\pi"] \arrow[d,"\widetilde p"'] & X \arrow[d,"p"] \\
\widetilde B \arrow[r,"\widetilde\pi"'] & B.
\end{tikzcd}
\label{diag:seifert.commutative.diagram}
\end{equation}
\item Let $\mathcal F$ denote the partition of $\widetilde X$ whose elements are the connected components of
$\pi^{-1}\!\bigl(p^{-1}(b)\bigr)$ for $b\in B$. Then $\mathcal F$ is the fiber decomposition of $\widetilde p$,
and in particular every element of $\mathcal F$ is a circle.
\item Let $\tau\colon \widetilde X\to \widetilde X$ be the canonical involution of the double branched cover $\pi$.
Then $\tau$ preserves the partition $\mathcal F$ and hence descends to an involution
$\widetilde\tau\colon \widetilde B\to \widetilde B$ satisfying
$\widetilde\pi\circ \widetilde\tau=\widetilde\pi$.
Moreover, the fibers of $\widetilde\pi$ are exactly the $\widetilde\tau$-orbits. In particular,
$\widetilde\tau$ fixes precisely the points of $\widetilde B$ lying over the branching set $p(\sing(X))$.

\end{enumerate}
\end{pro}

\begin{proof}
We construct $\widetilde B$, $\widetilde p$, and $\widetilde\pi$ from the generalized Seifert fibration $p\colon X\to B$ and the double branched cover $\pi\colon \widetilde X\to X$, and then verify the stated properties.

\medskip

\subsubsection*{The lifted partition.} 
Let $\mathcal F$ be as in the statement of the proposition. 
We call elements of $\mathcal{F}$ \emph{fibers}.

\subsubsection*{The elements of $\mathcal{F}$ are circles.}
Let $b\in B$. By the definition of a generalized Seifert fibration, every fiber $p^{-1}(b)$ is either a circle (a $C$-fiber) or a closed interval (an $I$-fiber), and $I$-fibers occur precisely over $p(\sing(X))$.

Suppose that $b\notin p(\sing(X))$. Then $p^{-1}(b)$ is a circle contained in $X\setminus \sing(X)$. 
Choose a sufficiently small orbifold neighborhood $U\subset B$ of $b$ disjoint from $p(\sing(X))$.
By the local description of generalized Seifert fibrations, the restriction
$p|_{p^{-1}(U)}\colon p^{-1}(U)\to U$ is fiber-preserving homeomorphic to a (classical) Seifert fibration of a solid torus (see Section~\ref{ss:generalized.seifert.fiber.spaces}).
Since $p^{-1}(U)\subset X\setminus\sing(X)$, the restriction $\pi|_{\pi^{-1}(p^{-1}(U))}\colon \pi^{-1}(p^{-1}(U))\to p^{-1}(U)$ is the orientation double cover of $p^{-1}(U)$. 
Thus, since $p^{-1}(U)$ is orientable, $\pi|_{\pi^{-1}(p^{-1}(U))}$ is the trivial double cover, so $\pi^{-1}(p^{-1}(U))$ consists of two connected components, each mapped homeomorphically onto $p^{-1}(U)$ by $\pi$.
In particular, $\pi^{-1}(p^{-1}(b))$ consists of exactly two connected components, each mapped homeomorphically onto $p^{-1}(b)$ by $\pi$. Since $p^{-1}(b)$ is a circle, each component is a circle. 
Therefore, every connected component of $\pi^{-1}(p^{-1}(b))$, that is, every element of $\mathcal F$ lying over $b$, is a circle.

Suppose now that $b\in p(\sing(X))$.
Then $p^{-1}(b)$ is an $I$-fiber, and there exists an orbifold neighborhood $U\subset B$ of $b$
such that $p^{-1}(U)$ is fiber-preserving homeomorphic to $\bpt$ with respect to the standard fibration $q\colon \bpt \to K_1(S^1(1/2))$ sending the $I$-fiber $q^{-1}(o)$ to the vertex $o$.
By the definition of the canonical orientable double branched cover, $\pi\colon \widetilde{X}\to X$ is a twofold branched cover with branching set $\sing(X)$. 
Under the fiber-preserving homeomorphism $p^{-1}(U)\cong \bpt$, which identifies $\sing(X)\cap p^{-1}(U)$ with 
$\sing(\bpt)$, the restriction $\pi|_{\pi^{-1}(p^{-1}(U))}\colon \pi^{-1}(p^{-1}(U))\to p^{-1}(U) $ is isomorphic to the model double branched cover $S^1\times D^2\to \bpt$ described in Lemma~\ref{L:double.branched.cover} (via a $\mathbb{Z}_2$-equivariant homeomorphism of total spaces which, after identifying $p^{-1}(U)$ with $\bpt$, commutes with the covering maps). 
Under this identification, the fiber $p^{-1}(b)$ corresponds to the $I$-fiber $q^{-1}(o)$.
In the model double branched cover $S^1\times D^2\to \bpt$ from Lemma~\ref{L:double.branched.cover}, where $q\colon \bpt\to K_1(S^1(1/2))$ is induced by $(e^{i\theta},x)\mapsto x$, the preimage of the $I$-fiber $q^{-1}(o)$ is the core circle $S^1\times \{0\}$ in $S^1\times D^2$. Hence, $\pi^{-1}(p^{-1}(b))$ is a circle. 
Combined with the previous case, we conclude that for all $b\in B$, every element of $\mathcal F$ is a circle.

\subsubsection*{Local Seifert charts for $\mathcal{F}$}
Let $L\in\mathcal F$ be a fiber and pick $\widetilde x\in L$. Set $x:=\pi(\widetilde x)$ and $b:=p(x)$.
Choose an orbifold neighborhood $U\subset B$ of $b$ so that the restriction $p|_{p^{-1}(U)}\colon p^{-1}(U)\to U$ is fiber-preserving homeomorphic to one of the standard local models in the definition, either a (classical) Seifert fibration of a solid torus or the standard model $q\colon\bpt\to K_1(S^1(1/2))$.
We consider two cases, corresponding to whether $b\in p(\sing(X))$.

Suppose first that $b\notin p(\sing(X))$.
Then $p|_{p^{-1}(U)}\colon p^{-1}(U)\to U$ is fiber-preserving homeomorphic to a (classical) Seifert fibration of a solid torus. 
Moreover, the restriction $\pi|_{\pi^{-1}(p^{-1}(U))}\colon \pi^{-1}(p^{-1}(U))\to p^{-1}(U)$ is a regular double cover, since $U\cap p(\sing(X))=\varnothing$, and hence $p^{-1}(U)\subset X\setminus \sing(X)$.
As noted above, since $p^{-1}(U)$ is orientable, $\pi|_{\pi^{-1}(p^{-1}(U))}$ is the trivial double cover. 
Hence $\pi^{-1}(p^{-1}(U))$ has two connected components swapped by the deck involution $\tau\colon \widetilde{X}\to \widetilde{X}$, the canonical involution of the double branched cover $\pi$.
In particular, $\pi^{-1}(p^{-1}(b))$ has two connected components: two circles exchanged by $\tau$, each mapping homeomorphically onto $p^{-1}(b)$ under $\pi$.
Therefore, $\pi^{-1}(p^{-1}(U))$ is a $3$-manifold covered by (classical) Seifert charts, and the fibers of $\mathcal F$
restricted to $\pi^{-1}(p^{-1}(U))$ are exactly the fibers of the lifted Seifert charts.

Suppose now that $b\in p(\sing(X))$. 
Then $p|_{p^{-1}(U)}\colon p^{-1}(U)\to U$ is fiber-preserving homeomorphic to the standard generalized Seifert fibration $q\colon \bpt \to K_1(S^1(1/2))$ sending the $I$-fiber $q^{-1}(o)$ to the vertex $o$. Note that this homeomorphism identifies $\sing(X)\cap p^{-1}(U)$ with $\sing(\bpt)$. 
By the definition of the canonical orientable double branched cover $\pi$, the restriction $\pi|_{\pi^{-1}(p^{-1}(U))}\colon \pi^{-1}(p^{-1}(U))\to p^{-1}(U)$ is a double branched cover of $p^{-1}(U)$ with branching set $\sing(X)\cap p^{-1}(U)$. 
After identifying $p^{-1}(U)$ with $\bpt$ as above, Lemma~\ref{L:double.branched.cover} implies that the lifted partition by elements of $\mathcal F$ is locally given by a (classical) Seifert chart on a solid torus.
Under this identification, the elements of $\mathcal{F}$ in 
$\pi^{-1}(p^{-1}(U))$ coincide with the circle fibers of the standard Seifert fibration of the solid torus 
$S^1\times D^2$.
In particular, the quotient $\pi^{-1}(p^{-1}(U))/\mathcal{F}$ is a disk and the map induced by 
$\pi$ between the local bases (quotients by circle fibers) $\pi^{-1}(p^{-1}(U))/\mathcal{F}\to U$ is a twofold branched cover of $U\cong K_1(S^1(1/2))$, branched over $b\in U$ (corresponding to the vertex of $K_1(S^1(1/2))$).

In both cases, we obtain a neighborhood $V\subset \widetilde X$ of $\widetilde x$ saturated by fibers of $\mathcal F$,
together with a homeomorphism from $V$ to a (classical) Seifert chart modelled on
$(S^1\times D^2)/\mathbb Z_k$ for some $k\ge 1$, sending fibers of $\mathcal F$ to circle fibers.
Thus, $\mathcal F$ defines a Seifert fibration on the $3$-manifold $\widetilde X$.

\subsubsection*{Definition of $\widetilde B$ and $\widetilde p$}
Define the base space $\widetilde B := \widetilde X/\mathcal F$, endowed with the quotient topology, and let
$\widetilde p\colon \widetilde X \longrightarrow \widetilde B$ be the quotient map.
By our previous analysis, around every point of $\widetilde X$ there are Seifert charts in which $\widetilde p$ corresponds to the standard projection $(S^1\times D^2)/\mathbb{Z}_k\to D^2/\mathbb{Z}_k$. 
It follows that $\widetilde B$ admits a $2$-orbifold structure whose local charts are of the form $D^2/\mathbb{Z}_k$ (with $k\geq 1$ depending on the chart), and with respect to this structure $\widetilde p$ is a Seifert fibration.
The compactness of $\widetilde B$ follows from the compactness of $\widetilde X$.

\subsubsection*{Definition of $\widetilde\pi$ and commutativity of diagram \eqref{diag:seifert.commutative.diagram}}
Let us first verify that the continuous map $p\circ \pi\colon \widetilde X\to B$ is constant on the elements of $\mathcal F$.
By construction, each fiber $L\in\mathcal{F}$ is contained in $\pi^{-1}(p^{-1}(b))$ for some $b\in B$.
Hence
$p(\pi(\widetilde x))=b$ for all $\widetilde x\in L$.
Define $\widetilde\pi\colon \widetilde B \longrightarrow B$ by
\begin{equation}
\label{eq:definition.of.tilde.pi}
\widetilde\pi\bigl(\widetilde p(\widetilde x)\bigr) := p(\pi(\widetilde x)).
\end{equation}
Since $p\circ\pi$ is constant on elements of $\mathcal{F}$ (equivalently, on the fibers of $\widetilde{p}$), the map $\widetilde{\pi}$ is well defined. Indeed, if $\widetilde{p}(\widetilde{x}) = \widetilde{p}(\widetilde{y})$, then $\widetilde{x}$ and $\widetilde{y}$ lie in the same fiber $L\in \mathcal{F}$. Hence $p(\pi(\widetilde{x}))=p(\pi(\widetilde{y}))$.
This gives the commutative diagram in the statement.


\subsubsection*{Induced involution on the base and fibers of $\widetilde\pi$.}
Let $\tau\colon \widetilde X\to \widetilde X$ denote the canonical involution of the double branched cover $\pi$. 
Hence, $\pi\circ \tau=\pi$ and $\tau^2=\mathrm{id}_{\widetilde X}$, and the fixed point set of $\tau$ is the branching set of $\pi$ (see the discussion of the canonical double branched cover in Section~\ref{ss:alexandrov.3-spaces}).
Since $\pi\circ\tau=\pi$, for every $b\in B$, we have
\[
\tau\bigl(\pi^{-1}(p^{-1}(b))\bigr)=\pi^{-1}(p^{-1}(b)).
\]
Since $\tau$ is a homeomorphism, it permutes the connected components of $\pi^{-1}(p^{-1}(b))$.
Equivalently, $\tau$ preserves the partition $\mathcal{F}$. 
Therefore, $\tau$ induces a continuous map $\widetilde\tau\colon \widetilde B\to \widetilde B$ given by 
\[
\widetilde\tau\bigl(\widetilde p(\widetilde x)\bigr):=\widetilde p\bigl(\tau(\widetilde x)\bigr),
\qquad \widetilde x\in \widetilde X.
\]
Note that $\widetilde{\tau}$ is well defined because $\tau$ sends $\mathcal{F}$-fibers to $\mathcal{F}$-fibers. 
Moreover, $\widetilde\tau$ is an involution since $\widetilde\tau^2(\widetilde p(\widetilde x))=\widetilde p(\tau^2(\widetilde x))=\widetilde p(\widetilde x)$.
In particular, $\widetilde{\tau}$ is a homeomorphism.
Finally, using the definition of $\widetilde\pi$ in \eqref{eq:definition.of.tilde.pi} and $\pi\circ\tau=\pi$, we obtain
\[
\widetilde\pi\bigl(\widetilde\tau(\widetilde p(\widetilde x))\bigr)
=\widetilde\pi\bigl(\widetilde p(\tau(\widetilde x))\bigr)
=p\bigl(\pi(\tau(\widetilde x))\bigr)
=p\bigl(\pi(\widetilde x)\bigr)
=\widetilde\pi\bigl(\widetilde p(\widetilde x)\bigr).
\]
Hence, $\widetilde\pi\circ \widetilde\tau=\widetilde\pi$.

We now identify the fibers of $\widetilde\pi$ with $\widetilde\tau$-orbits.
Fix $b\in B$ and consider the set $\widetilde\pi^{-1}(b)\subset \widetilde B$.
By definition,
\[
\widetilde\pi^{-1}(b)
=\Bigl\{\widetilde p(L)\ :\ L \text{ is a connected component of } \pi^{-1}(p^{-1}(b))\Bigr\}.
\]
If $b\notin p(\sing(X))$, then $p^{-1}(b)$ is a $C$-fiber contained in $X\setminus \sing(X)$. 
Hence, $\pi^{-1}(p^{-1}(b))$ is the disjoint union of two circles, and the involution $\tau$ exchanges these two components.
Therefore, $\widetilde\pi^{-1}(b)$ consists of two points $\{u,\widetilde\tau(u)\}$ with $\widetilde\tau(u)\neq u$.
If instead $b\in p(\sing(X))$, then $p^{-1}(b)$ is an $I$-fiber and, by Lemma~\ref{L:double.branched.cover},
$\pi^{-1}(p^{-1}(b))$ is connected. 
Hence, $\widetilde\pi^{-1}(b)$ consists of a single point $u$, which is fixed by $\widetilde\tau$.
In particular, for all $u,v\in \widetilde B$, we have $\widetilde\pi(u)=\widetilde\pi(v)$ if and only if
$v\in\{u,\widetilde\tau(u)\}$, i.e., the fibers of $\widetilde\pi$ are exactly the $\widetilde\tau$-orbits, with fixed points precisely over the branch points of $\widetilde{\pi}$, that is, over $p(\sing(X))$.

\subsubsection*{The map $\widetilde\pi\colon \widetilde{B}\to B$ is a twofold branched cover with branching set $p(\sing(X))$}
Let $b\in B$.
Suppose first that $b\notin p(\sing(X))$.
Choose a sufficiently small orbifold neighborhood $U\subset B$ of $b$ disjoint from $p(\sing(X))$, so $p^{-1}(U)\subset X\setminus \sing(X)$. 
Then the restriction $\pi|_{\pi^{-1}(p^{-1}(U))}\colon \pi^{-1}(p^{-1}(U))\to p^{-1}(U)$ is a regular twofold cover. By our previous analysis, $\pi\colon \widetilde{X}\to X$ maps each $\mathcal{F}$-fiber in $\pi^{-1}(p^{-1}(U))$ onto a $p$-fiber in $p^{-1}(U)$.
Passing to the quotients by the circle fibers of the Seifert fibrations $p|_{p^{-1}(U)}\colon p^{-1}(U)\to U$ and $\widetilde{p}|_{\widetilde{p}^{-1}(\widetilde{\pi}^{-1}(U))}\colon \widetilde{p}^{-1}(\widetilde{\pi}^{-1}(U))\to \widetilde{p}^{-1}(\widetilde{\pi}^{-1}(U))/\mathcal{F} = \widetilde{\pi}^{-1}(U)$
yields that $\widetilde\pi|_{\widetilde\pi^{-1}(U)}\colon\widetilde\pi^{-1}(U)\to U$ is an unbranched twofold covering of orbifolds.

Suppose now that $b\in p(\sing(X))$. 
Choose a sufficiently small neighborhood $U$ of $b$ so that the restriction $p|_{p^{-1}(U)}\colon p^{-1}(U)\to U$ is fiber-preserving homeomorphic to the standard generalized Seifert fibration $q\colon \bpt \to K_1(S^1(1/2))$ sending the $I$-fiber $q^{-1}(o)$ to the vertex $o$.
By our previous analysis, the induced map on the local bases is a twofold branched cover  $D^2\to K_1(S^1(1/2))$,
branched at its vertex. Hence, $\widetilde\pi$ is branched over $b$.

By the first case, $\widetilde{\pi}\colon \widetilde{B}\to B$ is unbranched over $B\setminus p(\sing(X))$. 
The second case shows $\widetilde{\pi}$ is branched over each $b\in p(\sing(X))$. Hence, the branching set of $\widetilde{\pi}$ is $p(\sing(X))$.
\medskip

Thus, we have shown that $\widetilde p$ is a Seifert fibration $\widetilde X\to \widetilde B$,
that $\widetilde\pi\colon \widetilde B\to B$ is a twofold branched cover with branching set $p(\sing(X))$,
and that $\mathcal F$ is precisely the fiber decomposition of $\widetilde p$. 
This verifies item (3) by the construction of $\widetilde B$ and $\widetilde p$, item (2) by \eqref{eq:definition.of.tilde.pi}, item (4) by the induced involution argument, and item (1) by the final local analysis.
\end{proof}


\subsubsection*{The lifted Seifert invariants.} We now determine the Seifert invariants of the Seifert fibration $\tilde{p}\colon \widetilde{X}\to \widetilde{B}$.
Let us denote the symbolic invariants associated to $X=\seif(B)$ as in \eqref{list-of-invariants} and those of $\tilde{X}=\seif(\tilde{B})$ by
\begin{equation}
\left\{ \tilde{b},\tilde{\varepsilon},\tilde{g}, \tilde{\iota}, \left\{ (\tilde{\alpha}_i, \tilde{\beta}_i) \right\}_{i=1}^{\tilde{n}} \right\}.
\end{equation}
It is then immediately clear that, since $\sing(\tilde{X})=\emptyset$, then $\tilde{\iota}=0$. Near the exceptional $C$-orbits, the canonical involution of $\tilde{X}$ acts freely and, therefore, the preimage in $\tilde{X}$ of each exceptional orbit on $X$ consists of two copies of itself. Thus, $\tilde{n}=2n$ and we can enumerate the Seifert invariants of $\tilde{X}$, for example, as $(\tilde{\alpha}_{2k},\tilde{\beta}_{2k})=(\alpha_k, \beta_k)$ and $(\tilde{\alpha}_{2k-1},\tilde{\beta}_{2k-1})=(\alpha_k, \beta_k)$ for each $k=1,2,\ldots,n$. Because the local actions of $\mathbb{Z}_2$ on $\tilde{B}$ that give rise to $B$ are orientation preserving, it follows that $\tilde{B}$ is orientable if and only if $B$ is orientable. Therefore, $\tilde{\varepsilon} = \varepsilon$.

We now compute $\chi(\tilde B)$. Let $S\subset B$ be a disjoint union of $\iota$ small open $2$–disks, one around each point of $I$–fiber type, and let $B^\circ=B\setminus S$. The restriction of $\tilde\pi$ to 
$\tilde{B}^\circ=\tilde\pi^{-1}(B^\circ)$ is a regular double covering
\[
\tilde\pi\colon \tilde{B}^\circ \longrightarrow\ B^\circ.
\]
Hence,
\[
\chi(\tilde{B}^\circ)=2\chi(B^\circ)=2\bigl(\chi(B)-\iota\bigr).
\]
To recover $\chi(\tilde{B})$ from $\chi(\tilde{B}^\circ)$ we glue back the $\iota$ disks around the branch points. Each removed base disk has a single preimage disk in $\tilde B$. Therefore
\[
\chi(\tilde{B}) = \chi(\tilde{B}^\circ)+\iota=2\chi(B)-\iota.
\]
When $B$ is closed and orientable, $\chi(B)=2-2g$.
 By the previous paragraph, $\tilde{B}$ is also  orientable. 
Hence $\chi(\tilde{B})=2-2\tilde g$. Therefore, $\tilde{g}=2g+\frac{\iota}{2}-1$ and, in particular, $\iota$ must be even.

To compute $\tilde b$, we proceed as follows. 
Let $B_0\subset B$ be given by removing small disjoint open disks around the $\iota>0$ points of $I$-fiber type.
Set 
\[
X_0=p^{-1}(B_0) = X\setminus\bigcup_{j=1}^\iota \mathrm{int}B(\mathrm{pt})_j.
\]
Let $Y$ be the closed Seifert $3$-manifold obtained from $X_0$ by gluing, along each boundary Klein bottle, a twisted $F$-block (i.e., a solid Klein bottle) fiber-preservingly. Thus, we obtain $Y$ from $X$ by replacing each $B(\mathrm{pt})$ with a twisted $F$--block.

Consider the canonical double branched cover $\pi\colon \tilde{X}\to X$. Over $X_0$, the cover restricts to a regular double cover $\pi\colon \tilde X_0\to X_0$. Each boundary Klein  bottle $K^2\subset \partial X_0$ lifts to a boundary torus $T^2\subset \partial \tilde{X}_0$. The boundary coverings $T^2\to K^2$ are the same regardless of whether we attach a $\bpt$ or a twisted $F$-block: in both cases the boundary involution corresponds to the canonical orientable double covering $T^2\to K^2$. 

Fix boundary identifications compatible with the circle foliation on each $K^2 \subset \partial X_0$. The lifts of the boundary gluing maps used to build $X$ (attaching copies of $\bpt$) and $Y$ (attaching twisted $F$-blocks) coincide on $\partial \tilde{X}_0$, since in both cases the boundary covering is the canonical double cover $T^2\to K^2$. Consequently,  $\tilde X$ and the orientable double cover $\tilde Y$ are obtained from the $\tilde X_0$ by attaching the same collection of solid tori along the same boundary maps. Hence $\tilde X$ and $\tilde Y$ are fiber-preservingly homeomorphic. In particular, they determine the same Seifert fibration and invariants. 

Now we compute Euler numbers. We follow Scott's conventions (see \cite[Section 3]{Scott}). Since $\tilde{X}$ is a Seifert bundle over a closed orbifold with oriented total space, its (rational) Euler number is given by 
\begin{align}
\label{eq:euler.number.double.branched.cover}
e(\tilde{X})=-\left(\tilde{b}+\sum_{i=1}^{\tilde{n}}\frac{\tilde{\beta}_i}{\tilde{\alpha}_i}\right) =
-\left(\tilde{b}+2\sum_{i=1}^n\frac{\beta_i}{\alpha_i}\right).
\end{align}
The Euler number is defined to be zero if the total space of the Seifert bundle is non-orientable. Since $Y$ contains a twisted $F$-block, which is a solid Klein bottle and hence non-orientable, $Y$ is non-orientable. Thus, $e(Y)=0$. For a finite covering $\hat{M}\to M$ of Seifert bundles of degree $d$, let $l$ be the degree of the induced orbifold covering of bases and $m$ the degree with which a regular fiber of $\hat{M}$ covers a regular fiber of $M$, so $d=lm$. Then 
\[
e(\hat{M}) = \frac{l}{m} e(M)
\]
(see \cite[Theorem 3.6]{Scott} and compare with \cite[Theorem~3.3]{JankinsNeumann1983}).
Hence, since $\tilde Y$ is the orientable double cover of $Y$, we have
    \[
    e(\tilde Y) = 2e(Y) = 0
    \]
Since $\tilde X\cong \tilde Y$ as Seifert bundles, we get $e(\tilde X) = 0$, and it follows from \eqref{eq:euler.number.double.branched.cover}
that 
\[
\tilde{b} = -2\sum_{i=1}^n\frac{\beta_i}{\alpha_i}.
\]
\qed 
\medskip

We now present two examples of non-equivalent generalized Seifert fibrations on the non-manifold Alexandrov space  $\Susp(P^2)\#\Susp(P^2)$, the connected sum of two copies of the suspension of $P^2$. Here, the connected sum is defined by taking the connected sum in the manifold part of $\Susp(P^2)$ (cf.\ \cite{FrancoGalazLarranagaGuijarroHeil}).

\begin{exa}[The double of $B(\mathrm{pt})$]
\label{ex:first.fibration}
Let
\[
X = B(\mathrm{pt})\ \cup_{K^2}\ B(\mathrm{pt}),
\]
obtained by gluing two copies of $B(\mathrm{pt})$ along their common boundary $K^2$ via the identity map, which is fiber-preserving. Then $X$ is a closed generalized Seifert fiber space with base $B\cong S^2$ with exactly two points of $I$–fiber type and no exceptional $C$–fibers. In the notation of Theorem~\ref{T:Classification-generalized-Seifert}, the invariants of $X$ are
\[
\{\, b=0,\ \varepsilon = \text{(orientable)},\ g=0,\ \iota=2,\ \{(\alpha_i,\beta_i)\}_{i=1}^{0}=\emptyset \,\}.
\]

By Lemma~\ref{L:double.branched.cover}, the canonical double branched cover of each half is the projection
$S^1\times D^2 \to B(\mathrm{pt})$, and the boundary covering is $T^2\to K^2$. Gluing the two lifted boundaries via the identity map of $T^2$ (the fiber-preserving lift of the identity map on $K^2\cong \partial\bpt$) yields
\[
\tilde{X} \;=\; (S^1\times D^2)\ \cup_{T^2}\ (S^1\times D^2)\cong S^2\times S^1,
\]
with its standard Seifert fibration $\tilde p:S^2\times S^1\to S^2$. This illustrates Theorem~\ref{T:Double-cover-is-seifert}; the Seifert fibration of $\tilde
{X}$ commutes with the branched cover, and the invariants transform as stated:
\[
\tilde\iota=0,\qquad \qquad \tilde\varepsilon=\varepsilon,\qquad
\chi(\tilde B)=2\,\chi(B)-\iota=2\ \ (\text{hence }\tilde g=0).
\]
For the underlying topology, $X$, being the double of $\bpt$, is homeomorphic to $\Susp(P^2)\#\Susp(P^2)$ (see, for example, \cite[Lemma 5.1]{FrancoGalazLarranagaGuijarroHeil}). 
\end{exa}

\begin{exa}[Two $I$–fibers and one $C$–fiber $(2,1)$]
\label{ex:second.fibration}
Let $P\subset S^2$ be the $2$-sphere with three open disks removed and let $Z\to P$ be the $\mathrm{O}(2)$–bundle of $R$–fibers whose boundary has two twisted components (total-space boundary $K^2\sqcup K^2$) and one untwisted component (total-space boundary $T^2$). Define
\[
X = \bigl(Z\ \cup_{\,K^2\sqcup K^2}\ \bigl(B(\mathrm{pt})\sqcup B(\mathrm{pt})\bigr)\bigr)\ \cup_{\,T^2}\ V_{(2,1)},
\]
where $V_{(2,1)}$ denotes a Seifert solid torus of type $(\alpha,\beta)=(2,1)$. We glue the $K^2$ boundary components of $Z$ and $\bpt$ via the identity map, which is fiber-preserving with respect to the chosen fiber-preserving identifications of each boundary component.
We attach $V_{(2,1)}$ by a fiber-preserving homeomorphism
$\varphi\colon \partial V_{(2,1)} \to \partial Z$ such that $\varphi_{*}(\mu)=2q_Z+h_Z$ in $H_1(\partial Z;\mathbb{Z})$, where $\mu$ is the class of the meridian of $\partial V_{(2,1)}$ in $H_1(\partial V_{(2,1)};\mathbb{Z})$,
$h_Z$ is the class of a regular fiber on $\partial Z$, and $q_Z$ is the class of a simple closed curve in $\partial Z$ transverse to the fibers such that $(q_Z,h_Z)$ is a $\mathbb{Z}$-basis of $H_1(\partial Z;\mathbb{Z})$ (see for example \cite[Lemma 6]{Seifert1933Eng}) or \cite[Theorem 1.10]{OrlikBook}).
Then $X$ is a closed generalized Seifert fiber space with base $B\cong S^2$ having two points of $I$–fiber type and one exceptional $C$–fiber of type $(2,1)$. In the notation of Theorem~\ref{T:Classification-generalized-Seifert}, the invariants of $X$ are
\[
\{\, b=0,\ \varepsilon=\text{(orientable)},\ g=0,\ \iota=2,\ \{(2,1)\}\,\}.
\]

By Lemma~\ref{L:double.branched.cover}, the canonical double branched cover of each $B(\mathrm{pt})$ is $S^1\times D^2$, with boundary covering $T^2\to K^2$. The Seifert solid torus $V_{(2,1)}$ does not intersect the branching set, so its preimage is the disjoint union of two Seifert solid tori of the same type $(2,1)$. 
Gluing the lifted pieces we obtain a Seifert fibration of the canonical double branched cover $\tilde X$ with base $\tilde{B}\cong S^2$
with two exceptional fibers $(2,1)$ and $\tilde b=-1$. Hence, $e(\tilde X)=-(\tilde b+\tfrac12+\tfrac12)=0$. In particular,
\[
\tilde X \cong M\bigl(S^2;\ \tilde b=-1;\ (2,1),(2,1)\bigr)\cong S^2\times S^1,
\]
as recorded in \cite[Section 4]{Scott}. This illustrates Theorem~\ref{T:Double-cover-is-seifert}: in $\tilde{X}$, the list of exceptional pairs is doubled and the Seifert fibration commutes with the branched covering. For the underlying topology, $X$ is the quotient of  $S^2\times S^1$ by an orientation-reversing involution with exactly four isolated fixed points.
This involution is unique up to conjugacy (see \cite[Theorem A]{Tollefson1979}) and its quotient is homeomorphic to $\Susp(P^2)\# \Susp(P^2)$ (see also \cite{GalazGuijarro3dim} and \cite[Lemma~5.1]{FrancoGalazLarranagaGuijarroHeil}). The double branched cover $\pi\colon \tilde{X}\to X$ has exactly four branching points, corresponding to the four fixed points of the involution on $X\cong S^2\times S^1$, and each one of the four branching points corresponds to an endpoint of the two $I$-fibers in $X\cong\Susp(P^2)\#\Susp(P^2)$.
\end{exa}


\bibliographystyle{amsplain}
	\bibliography{main}
 
\end{document}